\documentclass[10pt,a4paper,final]{article}
\title{Channel linear Weingarten surfaces}
\date{13 July 2015}
\author{U. Hertrich-Jeromin, K. Mundilova, E.-H. Tjaden}

\usepackage{amsmath}
\usepackage{amssymb}
\usepackage{amsthm}
\usepackage{graphicx}
\usepackage{subfig}

\newcommand{\R}{\mathbb{R}}

\newcommand{\cn}{\operatorname{cn}}
\newcommand{\dn}{\operatorname{dn}}
\newcommand{\sn}{\operatorname{sn}}

\newcommand{\am}{\operatorname{am}} 
\newcommand{\III}{{\rm I\!I\!I}}
\newcommand{\II}{{\rm I\!I}}
\newcommand{\I}{{\rm I}}

\newtheorem{theorem}{Theorem}[section]
\newtheorem{lemma}{Lemma}[section]
\newtheorem{corollary}{Corollary}[section]

\begin{document}
\maketitle

This note is based on the second author's Bachelor thesis,
the purpose of which was to understand the classification
of cyclic linear Weingarten surfaces from \cite{lo08}.
In particular, we obtained a very simple proof for Cor~3.6
of \cite{lo08}, using explicit parametrizations in terms of
Jacobi elliptic functions based on \cite{Tjaden91}.
This method will also be applicable to the results
of \cite{Lopez08-2}.
Moreover, in our attempt to derive a more conceptual proof
for the key step \cite[Thm~2.1]{lo08} in the classification
of cyclic linear Weingarten surfaces of \cite[Thm~1.1]{lo08},
we derived the classification of channel linear Weingarten
surfaces Thm~\ref{thmchannel} below
---
that, in fact, had already been obtained in \cite{Havlicek53},
using different methods based on \cite{Havlicek49}.
As a consequence of our Thm~\ref{thmchannel} we obtain a partial
but rather explicit classification,
of channel linear Weingarten surfaces as special cyclic linear
Weingarten surfaces, in Sect~3 of this note.

Though this paper does not contain substantial new results,
but merely employs various well known methods and results in
order to elucidate the main classification result of \cite{lo08},
we feel that its publication may serve the mathematical community
by recalling these methods and by demonstrating how they beautifully
serve to classify channel linear Weingarten surfaces;
in fact, basic scholarly work suggests that the classical and
rather direct methods we employ have fallen into oblivion
and require revivification.

\section{Parallel linear Weingarten surfaces}
It is well known that the parallel surfaces $x^t=x+t\,n$ of
a linear Weingarten surface $x$ in Euclidean space are linear
Weingarten:

\begin{lemma}\label{plw}
If the Gauss and mean curvatures $K$ and $H$, respectively,
of a surface $x:\Sigma^2\to\R^3$ satisfy a linear Weingarten
condition
\begin{equation}\label{lwc}
  0 = a\,K + 2b\,H + c
\end{equation}
then the Gauss and mean curvatures of its parallel surfaces
$x^t=x+t\,n$ satisfy
\begin{equation}\label{parallel}
  0 = (a + 2tb + t^2c)\,K^t + 2(b+tc)\,H^t +c.
\end{equation}
\end{lemma}

Note that surfaces of constant mean or constant Gauss curvature
($a=0$ or $b=0$, respectively)
yield examples of linear Weingarten surfaces, with parallel
linear Weingarten surfaces that do generally not have
a constant curvature.

A simple and fairly efficient proof of Lemma~\ref{plw}
is obtained by employing the Cayley-Hamilton identity:
for the fundamental forms of $x^t=x+t\,n$,
$$
  \III^t = \III, \enspace
  \II^t = \II - t\,\III \enspace{\rm and}\enspace
  \I^t = \I - 2t\,\II +  t^2\III,
$$
the Cayley-Hamilton identity $0=\III-2H\,\II+K\,\I$
of the initial surface $t=0$ yields
$$
  0 = (1-2Ht+Kt^2)\,\III^t - 2 (H-tK)\,\II^t + K\,I^t.
$$
Hence $K^t=\frac{K}{1-2tH+t^2K}$ and
$H^t=\frac{H-tK}{1-2tH+t^2K}$ so that the claim
follows from the linear Weingarten condition for the initial
surface,
$$
    (a + 2tb + t^2c)\,K
  + 2(b+tc)\,(H-t\,K)
  + c\,(1 - 2Ht + Kt^2)
  = 0.
$$

Note that the linear Weingarten condition \eqref{lwc}
factorizes when its discriminant vanishes,
$$
  0 = \Delta := -ac+b^2
   \enspace\Rightarrow\enspace
  0 = (a\kappa_1+b)(a\kappa_2+b),
$$
so that one of the principal curvatures $\kappa_i$ is constant,
$\kappa_i\equiv-\frac{b}{a}$,
cf \cite[\S2.V]{bo60}.
In this case the linear Weingarten surface $x$ is either
developable, if $b=0$, or it is a tube of radius $\frac{a}{b}$
about its focal curve $x-\frac{a}{b}\,n$;
in either case, we will refer to the linear Weingarten surface
$x$ as \emph{tubular}.
Note that we will also consider totally umbilic surfaces as
\emph{tubular linear Weingarten surfaces}.

If, on the other hand, $x$ is a non-tubular linear Weingarten
surface, $\Delta\neq0$,
then \eqref{parallel} shows that its parallel family
$(x+t\,n)_{t\in\R}$ contains a constant curvature
representative.

In particular,
if $c=0$ then $b\neq 0$ and $a+2tb=0$ yields $H^t=0$,
 that is, a minimal surface;
if $c\neq0$ then $b+tc=0$ yields a surface of constant
 Gauss curvature $K^t=\frac{c^2}{\Delta}$;
if, moreover, $\Delta>0$ then we recover Bonnet's
theorem on parallel surfaces of constant curvature from
\eqref{parallel} with $a+2tb+t^2c=0$,
that is, with $t=\frac{-b\pm\sqrt{\Delta}}{c}$,
cf \cite[\S2.7.4]{imdg}.

Thus we obtain three classes of parallel families of
non-tubular linear Weingarten surfaces:

\begin{corollary}\label{thmparallel}
Let $x:\Sigma^2\to\R^3$ be a non-tubular linear Weingarten
surface,
$$
  0 = a\,K + 2b\,H + c
   \enspace{with}\enspace
  \Delta := b^2-ac \neq 0.
$$
If $c=0$ then $x$ is parallel to a minimal surface;
if $c\neq0$ then $x$ is parallel to a surface of
 constant Gauss curvature $K\neq 0$,
which has two parallel constant mean curvature surfaces
 when $K>0$, i.e., when $\Delta>0$.
\end{corollary}

\section{Channel linear Weingarten surfaces}
Next we demonstrate that a linear Weingarten surface with
a family of circular curvature lines,
that is,
a \emph{channel} linear Weingarten surface,%
\footnote{A channel surface can be characterized as the envelope
 of a $1$-parameter family of spheres or,
 equivalently (in conformal geometry),
 as a surface with a family of circular curvature lines.
 \label{charchannel}}
is tubular or a surface of revolution,
cf \cite[Thm~1]{Havlicek53} or \cite[Thm~1.1]{lo08}.

Following \cite{lo08} we parametrize a surface foliated
by circles:
if $u\mapsto z(u)$ denotes the centre curve of the foliating
circles and $r=r(u)$ their radii%
\footnote{Here we exclude surfaces where the family
 of circular curvature lines contains lines,
 such as on developable surfaces.}
then
\begin{equation}\label{channel}
  (u,v) \mapsto x(u,v) = z(u) + r(u)\,\nu(u,v),
\end{equation}
where we choose standard parametrizations $v\mapsto\nu(u,v)$
of the unit circles, such that $\nu_v\perp\nu$ and,
with unit normals $\tau=\tau(u)$ of the planes of the foliating
circles,
we obtain an orthonormal frame field $(\tau,\nu,\nu_v)$
for the surface.
Assuming that $u$ parametrizes the circle planes by means
of a unit speed orthogonal trajectory,%
\footnote{Such an orthogonal trajectory does not always exist,
 as the example of a surface of revolution with closed profile
 curve shows.}
$\tau$ becomes the tangent field of this orthogonal trajectory.
To facilitate our analysis we require that $u\mapsto\nu(u,0)$,
hence $u\mapsto\nu(u,v)$ for every $v$,
be parallel along this orthogonal trajectory,
that is,
$$
  \nu_u = -\kappa\tau
   \enspace{\rm and}\enspace
  \tau' = \kappa\nu + \kappa_v\nu_v,
$$
where $\kappa$ and $\kappa_v$ are the curvatures of the
trajectory with respect to the (parallel) normal fields
$\nu$ and $\nu_v$, for fixed $v$.
Note that $\tau'_v=0$, hence $\kappa+\kappa_{vv}=0$.

As $\nu_v\parallel x_v$ is tangential, the Gauss map $n$
of the channel surface $x$ is given by
$$
  n = \tau\cos\alpha + \nu\sin\alpha,
   \enspace{\rm where}\enspace
  \alpha=\alpha(u)
$$
denotes the intersection angle of the surface with the planes
of the foliating circles
---
which is constant along the intersection lines $v\mapsto x(u,v)$
by Joachimsthal's theorem, since these lines are curvature lines
on the surface.
Note that this also proves one direction of the equivalence
noted in Footnote~\ref{charchannel}.
As a consequence, since $x_u\perp n$,
$$
  x_u = w\,(-\tau\sin\alpha+\nu\cos\alpha)
      + w_v\nu_v\cos\alpha,
   \enspace{\rm where}\enspace
  (u,v)\mapsto w(u,v)
$$
denotes a suitable function,
and the shape of the $\nu_v$-coefficient is obtained
by using that $0=z_{uv}=(x-r\nu)_{uv}$,
as are
\begin{equation}\label{coeff}
  r' = (w_{vv}+w)\cos\alpha
   \enspace{\rm and}\enspace
  \kappa_v = \tfrac{1}{r}w_v\sin\alpha,
   \enspace{\rm hence}\enspace
  \kappa = -\tfrac{1}{r}w_{vv}\sin\alpha
\end{equation}
since $\kappa_{vv}+\kappa=0$.
With $x_v=r\nu_v$ and
$$
  n_u = (\alpha'+\kappa)(-\tau\sin\alpha+\nu\cos\alpha)
      + \kappa_v\nu_v\cos\alpha
   \enspace{\rm and}\enspace
  n_v = \nu_v\sin\alpha
$$
it is now straightforward to determine
the Gauss and mean curvatures of $x$:
$$
  K = \tfrac{\alpha'+\kappa}{w} \tfrac{\sin\alpha}{r}
   \enspace{\rm and}\enspace
  H = -\tfrac{1}{2}(
    \tfrac{\alpha'+\kappa}{w} + \tfrac{\sin\alpha}{r}
    ).
$$
Hence the linear Weingarten condition $0=a\,K+2b\,H+c$
for our channel surface
$x$ reads
\begin{equation}\label{clw}
  0 = (\alpha'+\kappa)(a\sin\alpha-br) - w\,(b\sin\alpha-cr).
\end{equation}

Differentiating the linear Weingarten equation \eqref{clw}
twice and using \eqref{coeff} we deduce
$$
  0 = w_{vv}(a\sin^2\alpha-2br\sin\alpha+cr^2)
$$
as a necessary condition for a channel surface \eqref{channel}
to be linear Weingarten:
this yields two cases to analyze.

If $a\sin^2\alpha-2br\sin\alpha+cr^2\equiv 0$ on an open
interval then the channel surface $x$ has a constant
principal curvature,
$\tfrac{1}{r}\sin\alpha\equiv const$,
that is, $x$ is tubular on that interval.

If $a\sin^2\alpha-2br\sin\alpha+cr^2$ has at most isolated
zeroes then $w_{vv}\equiv 0$;
consequently, by \eqref{coeff},
$$
  \kappa \equiv 0,
   \enspace{\rm hence}\enspace
  w_v \equiv 0,
   \enspace{\rm and}\enspace
  r' = w\cos\alpha
$$
which show that the orthogonal trajectory of the circle planes
is a straight line, so that the circle planes are parallel,
and, with \eqref{channel} and $\kappa\equiv 0$,
$$
  z' = (x-r\nu)_u = -w\sin\alpha\,\tau
$$
so that the curve of circle centres is a straight line and,
in particular, an orthogonal trajectory of the circle planes.
Thus the channel surface $x$ is a surface of revolution
with its centre curve $z$ taking values on the axis,
cf \cite[\S2.III]{bo60} and \cite[Thm~1]{Havlicek53}:

\begin{theorem}\label{thmchannel}
A non-tubular linear Weingarten surface with circular curvature
lines is a surface of revolution.
\end{theorem}

\section{Linear Weingarten surfaces of revolution}
With Cor~\ref{thmparallel}, Thm~\ref{thmchannel}
yields a classification result for channel linear Weingarten
surfaces:

\begin{corollary}
A channel linear Weingarten surface is either tubular
(tube about a space curve, developable or totally umbilic)
or parallel surface of
 a minimal ($H=0$) or
 a constant Gauss curvature $K\neq 0$
surface of revolution.
\end{corollary}

It is well known that the catenoid is the only minimal surface
of revolution in Euclidean $3$-space, cf \cite[\S2.VIII]{bo60},
hence the following classification of surfaces of revolution with
constant Gauss curvature completes our classification
of channel linear Weingarten surfaces,
cf \cite[Sect 6.7]{Tjaden91}:

\begin{theorem}\label{thmrevolution}
Every constant Gauss curvature $K\neq 0$ surface of revolution
is homothetic to
$$
  (s,\vartheta) \mapsto x(s,\vartheta) :
  = (r(s)\cos\vartheta, r(s)\sin\vartheta, h(s)),
$$
with one of the following profile curves:
$$\begin{matrix}\hfill
 r(s)=p\cn_p(s), & \begin{cases}
  h(s)=(E_p\circ\am_p)(s) & (K=+1), \cr
  h(s)=(E_p\circ\am_p)(s)-s & (K=-1); \cr
  \end{cases} \hfill\cr\cr\hfill
 r(s)=\tfrac{1}{\cosh(s)}, & \begin{cases}
  h(s)=\tanh(s) & (K=+1,\text{ sphere}), \cr
  h(s)=\tanh(s)-s & (K=-1,\text{ pseudosphere}); \cr
  \end{cases} \hfill\cr\cr\hfill
 r(s)=\tfrac{1}{p}\dn_p(\tfrac{s}{p}), & \begin{cases}
  h(s)=\tfrac{1}{p}(E_p\circ\am_p)(\tfrac{s}{p})-\tfrac{1-p^2}{p^2}s &
   (K=+1), \cr
  h(s)=\tfrac{1}{p}(E_p\circ\am_p)(\tfrac{s}{p})-\tfrac{1}{p^2}s &
   (K=-1). \cr
  \end{cases} \hfill\cr
\end{matrix}$$
Here $\sn_p=\sin\circ\am_p$, $\cn_p=\cos\circ\am_p$ and $\dn_p$
denote Jacobi elliptic functions,
with the Jacobi amplitude function $\am_p$, and
$E_p(\phi)=\int_0^\phi\sqrt{1-p^2\sin^2\varphi}\,d\varphi$
is the incomplete elliptic integral of the second kind,
with elliptic module $p\in(0,1)$.
\end{theorem}

It is straightforward to verify that the surfaces listed
in Thm~\ref{thmrevolution} do indeed have constant Gauss
curvatures $K=\pm 1$:
the Gauss curvature of a surface of revolution is given by
$$
  K = \tfrac{r'h''-r''h'}{\sqrt{r'^2+h'^2}^3}
      \tfrac{h'}{r\sqrt{r'^2+h'^2}}
  = -\tfrac{1}{2rr'}(\tfrac{r'^2}{r'^2+h'^2})'
$$
in terms of its profile curve $s\mapsto(r(s),h(s))$;
hence the surface has constant Gauss curvature if and only if
there is a constant $C\in\R$ so that
$$
  r'^2 = (r'^2+h'^2) (C-Kr^2).
$$
Note that
$$
  0 \leq C-Kr^2 = \tfrac{r'^2}{r'^2+h'^2} \leq 1.
$$
Given a surface of revolution with constant Gauss curvature $K$
this differential equations for its profile curve can be further
simplified to an elliptic differential equation
\begin{equation}\label{DGLr}
  r'^2 = ((1-C)+Kr^2) (C-Kr^2)
\end{equation}
for the radius function $s\mapsto r(s)$ alone by adjusting
the parametrization of the profile curve so that
$$
  r'^2 + h'^2 = (1-C) + Kr^2
   \enspace{\rm as}\enspace
  (1-C) + Kr^2 \geq 0;
$$
hence, without loss of generality,
\begin{equation}\label{DGLh}
  h' = (1-C)+Kr^2.
\end{equation}

To analyze the solutions of the differential equations
\eqref{DGLr} and \eqref{DGLh} that yield the occurring
profile curves we distinguish the cases $K=\pm 1$.

If $K=+1$ then $0\leq r^2\leq C\leq 1+r^2$ and \eqref{DGLr} reads
$$
  y'^2 = \begin{cases}
   (1-y^2)(q^2+p^2y^2) \hfill&
    \text{with $p=\sqrt{C}$, $q=\sqrt{1-C}$
     and $y=\tfrac{r}{p}$ if $C<1$,}\hfill\cr
   y^2-1 \hfill&
    \text{with $y=\tfrac{1}{r}$ if $C=1$, and}\hfill\cr
   \tfrac{1}{p^2}(1-y^2)(y^2-q^2) \hfill&
    \text{with $p=\tfrac{1}{\sqrt{C}}$, $q=\sqrt{\tfrac{C-1}{C}}$
     and $y=pr$ if $C>1$;}\hfill\cr
  \end{cases}
$$
the solutions of these differential equations are,
up to parameter shift,
$$
  y(s)=\cn_p(s), \enspace
  y(s)=\cosh s \enspace{\rm and}\enspace
  y(s)=\dn_p(\tfrac{s}{p}),
$$
respectively.
Integration of \eqref{DGLh} then yields the respective three
profile curves, as claimed in Thm~\ref{thmrevolution},
up to translation or reflection.

If $K=-1$ then $0\leq r^2\leq 1-C\leq 1+r^2$ and swapping
the roles of $1-C$ and $C$ leads to the same elliptic
differential equations for $r$ as in the $K=+1$ cases;
up to a sign choice, \eqref{DGLh} then differs from the
same equation in the $K=+1$ case by an additive $1$.
Hence we obtain the profile curves as claimed
in Thm~\ref{thmrevolution},
up to translation and reflection again.

This completes the proof of Thm~\ref{thmrevolution}.

The explicit parametrizations of Thm~\ref{thmrevolution}
now enable us to verify the results of \cite{Lopez08-2}
and, in particular, to confirm Cor~3.6 of \cite{lo08}
in a simple way:

\begin{corollary}
There is a $1$-parameter family of complete, immersed,
rotational linear Weingarten surfaces of hyperbolic type,
that is, with $\Delta<0$.\par
These surfaces have a translational period.
\end{corollary}

\begin{figure}
\subfloat[rotational surface]{\includegraphics[width=5cm]{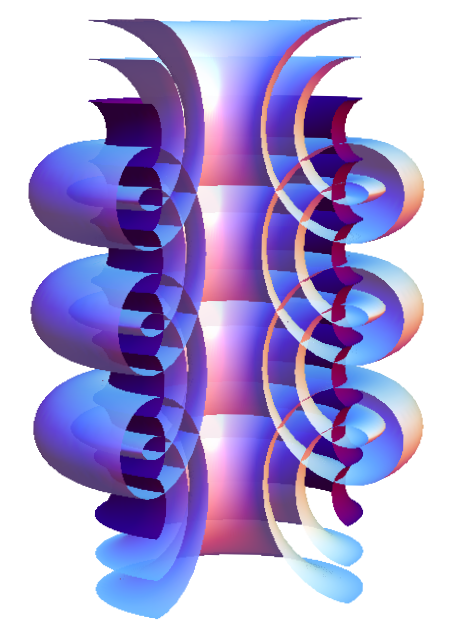}}
\subfloat[profile curves]{\includegraphics[width=7cm]{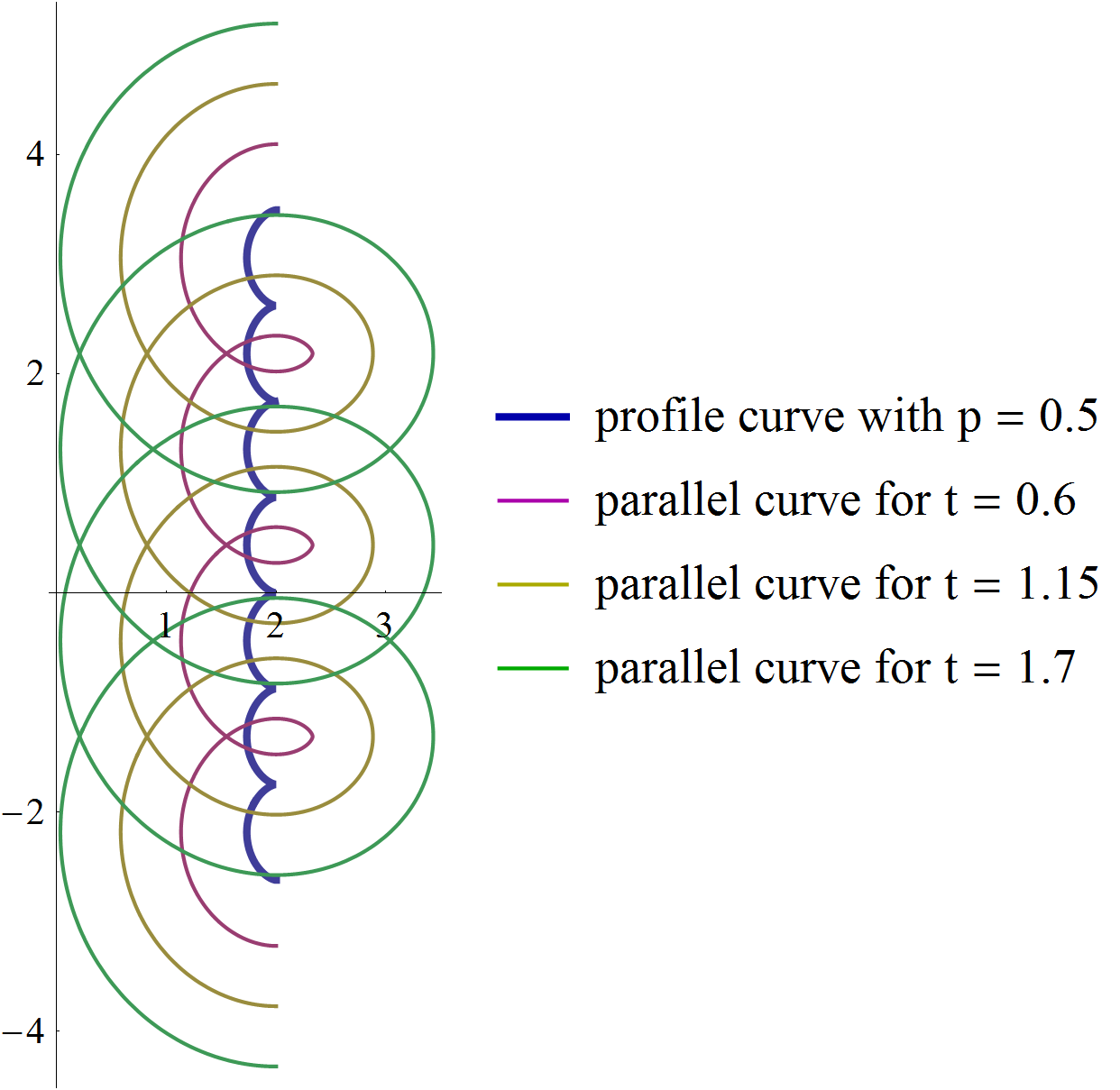}}
\caption{Parallel hyperbolic linear Weingarten surfaces of revolution}
\label{fig1}
\end{figure}

As the discriminant of the linear Weingarten equation is an
invariant of a parallel family of linear Weingarten surfaces,
$\Delta^t=\Delta$ by \eqref{parallel}, we consider the parallel
surfaces of a pseudospherical surface of revolution,
where $\Delta=-1$:
for
$$
  \xi = (r,h)
   \enspace{\rm with}\enspace
  r(s) = \tfrac{1}{p}\dn_p(\tfrac{s}{p})
   \enspace{\rm and}\enspace
  h(s) = \tfrac{1}{p}(E_p\circ\am_p)(\tfrac{s}{p})-\tfrac{1}{p^2}s
$$
we obtain
$\xi'(s)=-\sn_p(\tfrac{s}{p})\,\tau(s)$ with
$\tau(s)=(\cn_p,\sn_p)(\tfrac{s}{p})$;
hence the parallel curves
$$
  \xi^t = \xi + t\,\nu
   \enspace{\rm with}\enspace
  \nu(s) = (-\sn_p,\cn_p)(\tfrac{s}{p})
$$
hit the axis of rotation where
$$
  0 = r^t(s) = (\tfrac{1}{p}\dn_p-t\sn_p)(\tfrac{s}{p})
   \enspace\Leftrightarrow\enspace
  t = \tfrac{1}{p}\tfrac{\dn_p}{\sn_p}(\tfrac{s}{p})
   \in(-\infty,-\tfrac{q}{p}]\cup[\tfrac{q}{p},\infty)
$$
with the elliptic co-module $q=\sqrt{1-p^2}$,
and they develop singularities where
$$
  0 = (\xi^t)'(s)
  = -(\sn_p+\tfrac{t}{p}\dn_p)(\tfrac{s}{p})\,\tau(s)
   \enspace\Leftrightarrow\enspace
  t = -p\tfrac{\sn_p}{\dn_p}(\tfrac{s}{p})
   \in[-\tfrac{p}{q},\tfrac{p}{q}].
$$
Consequently, if the elliptic modulus is chosen so that
$$
  \tfrac{p}{q} < \tfrac{q}{p}
   \enspace\Leftrightarrow\enspace
  p^2 < q^2 = 1 - p^2
   \enspace\Leftrightarrow\enspace
  p < \tfrac{1}{\sqrt{2}},
   \enspace{\rm then}\enspace
  t \in (\tfrac{p}{q},\tfrac{q}{p}) \neq \emptyset
$$
yield a family of immersed parallel linear Weingarten surfaces
of revolution with profile curves $\xi^t$, cf Figure \ref{fig1}.
Clearly these surfaces are complete,
and they have a translational period since $r$ and $h$ are periodic.

\bibliographystyle{clw}

\end{document}